\documentclass[preprint,12pt]{elsarticle}

\usepackage{amssymb}
\usepackage{amsthm}

\usepackage{multirow}  			
\usepackage{amsmath}

\newtheorem{prop}{Proposition}
\newtheorem{lemma}{Lemma}
\newtheorem{theorem}{Theorem}

\newtheorem*{quest}{Question}
\newtheorem*{remark}{Remark}

\newdefinition{dfn}{Definition}

\newproof{pf}{Proof}

\def\FF{\mathbb F}

\def\PP{\mathbb P}

\def\QQ{\mathbb Q}

\def\ZZ{\mathbb Z}

\def\tensor{\otimes}
\def\bigtensor{\bigotimes}

\DeclareMathSymbol\circlearrowright{\mathrel}{AMSa}{"08} 
\DeclareMathSymbol\circlearrowleft{\mathrel}{AMSa}{"09} 
\DeclareMathSymbol\Vdash{\mathrel}{AMSa}{"0D} 
\DeclareMathSymbol\Vvdash{\mathrel}{AMSa}{"0E} 
\DeclareMathSymbol\vDash{\mathrel}{AMSa}{"0F} 
\DeclareMathSymbol\twoheadrightarrow{\mathrel}{AMSa}{"10} 
\DeclareMathSymbol\twoheadleftarrow{\mathrel}{AMSa}{"11} 

\def\onto{\twoheadrightarrow}
\def\into{\hookrightarrow}

\newcommand\twodcol[2]{ \ensuremath{ \left [
	\begin{array}{c}
		#1  \\
		#2
	\end{array}  
	\right ] } }
  
\newcommand\twodmatrix[4]{ \ensuremath{ \left( 
	\begin{array}{cc}
		#1 & #2  \\
		#3 & #4 
	\end{array}  
	\right) } }

\renewcommand{\mod}{\mathrm{ mod }}

\journal{Journal of Number Theory}

\begin{document}

\begin{frontmatter}



\title{On Serre's Conjecture over Imaginary Quadratic Fields}


\author{Rebecca Torrey\footnote{Present Address: Department of Mathematics, Amherst College Box 2239, Amherst MA 01002, rtorrey@mtholyoke.edu, (413) 542-2100 (Dept Office), (413) 542-2550 (FAX)}}

\address{Department of Mathematics, King's College London, Strand, London WC2R 2LS, England}

\begin{abstract}
We study an analog of Serre's conjecture over imaginary quadratic fields.  In particular, we ask whether the weight recipe of Buzzard, Diamond and Jarvis will hold in this setting.  Using a program written by the author, we provide computational evidence that this is in fact the case.  In order to justify the method used in the program, we prove that a modular symbols method will work for arbitrary weights over imaginary quadratic fields.  
\end{abstract}


\end{frontmatter}


\section {Introduction} \label{intro}

\medskip

Serre's conjecture \cite{serre2} states that any continuous, odd, irreducible representation
\[
	\rho : G_{\QQ} = \mathrm{Gal} (\overline{\QQ}/\QQ) \to \mathrm{GL} _2 (\overline{\FF}_\ell)
\]
is modular, i.e., arises from a modular form. A refinement to Serre's conjecture gives the minimal weight and level of this modular form and, by work of Ribet and others, we know that Serre's conjecture holds if and only if the refined version holds (assuming $\ell > 2$). Recently, Khare and Wintenberger \cite{khare_wintenberger_I, khare_wintenberger_II} have completed the proof of Serre's conjecture using ideas and results of Dieulefait \cite{dieulefait}, Kisin \cite{kisin1, kisin2}, Taylor \cite{taylor} and Wiles \cite{taylor_wiles, wiles}.  

It is natural to ask whether an analogous conjecture holds for representations of $G_K$ where $K$ is an arbitrary number field. Buzzard, Diamond and Jarvis \cite{bdj} recently formulated a version of the refined Serre's conjecture in the case in which $K$ is totally real and $\ell$ is unramified in $K$, where predicting the weights is much more complicated than for $K = \QQ$.  In this more general setting one needs a $d$-tuple to specify the weight (where $d = [K:\QQ]$), so it no longer makes sense to simply specify a minimal weight and level for the corresponding modular form.  A more general notion of weight is needed.  Also, the refined part of the conjecture takes the form of a recipe for {\it all} the weight combinations for modular forms giving rise to a particular representation.  This recipe depends only on the local behavior of the representation $\rho$ at primes above $\ell$.  The methods of Khare and Wintenberger break down in this situation, but there has been some progress, due to Gee \cite{gee1, gee2}, towards proving the equivalence of Serre's conjecture and the refined Serre's conjecture for totally real fields $K$.  If $K$ is not totally real, then the situation is even less well understood.  In this case, we do not even have a complete understanding of how to associate Galois representations to modular forms.  This situation is the focus of this paper. 

Some computational evidence is available for generalizations of Serre's conjecture to number fields.  Demb\'{e}l\'{e} \cite{ddr} has done computations of arbitrary weight mod $\ell$ modular forms over totally real fields $F$, providing evidence for the conjecture of Buzzard, Diamond and Jarvis.  For imaginary quadratic fields, Figueiredo \cite{figueiredo} provided some computational evidence for Serre's conjecture, but he worked only with weight two modular forms.  More recently, \c{S}eng\"{u}n \cite{sengun} proved non-existence of certain representations and has done some computations of arbitrary weight modular forms over imaginary quadratic fields. 

This paper poses the question of whether a conjecture analogous to that of Buzzard, Diamond and Jarvis will hold over imaginary quadratic fields.  In Section \ref{modforms}, it is shown that a modular symbols method akin to those used by Cremona et al. can be used to compute arbitrary weight modular forms over $\QQ(i)$ in this setting.   To provide evidence that a conjecture analogous to the BDJ conjecture will hold, examples of Galois representations are computed in Section \ref{galreps} and then code written by the author (employing the above-mentioned modular symbols method) is used to produce corresponding mod $\ell$ modular forms in the weights predicted by the BDJ recipe.

\section {Definitions} \label{defns}

\subsection {Serre Weights} \label{serreweights}

Buzzard, Diamond and Jarvis define Serre weights to account for the more complicated weight structure when $K$ is not $\QQ$.  They work with totally real fields, but both their definition of Serre weights and their recipe for the weight conjecture can be applied to general number fields.  

Let $K$ be a number field and denote by $\mathcal O$ its ring of integers.  Fix a prime $\ell$ that is unramified in $K$.  
\begin{dfn}
	We define a \emph{Serre weight} to be an irreducible $\bar \FF_\ell$-representation $V$ of $G = GL_2(\mathcal O / \ell \mathcal O )$.
\end{dfn}
Ash, Doud and Pollack \cite[p. 4]{ashDP} describe why this is a natural generalization of the weights in Serre's original conjecture.

We can describe Serre weights explicitly as follows: Set
\[
G = GL_2(\mathcal O / \ell \mathcal O ) \cong \prod_{\mathfrak p | \ell} GL_2(\mathcal O / \mathfrak p).
\]
For each prime $\mathfrak p$ of $K$ such that $\mathfrak p | \ell$, set $k_{\mathfrak p} = \mathcal O / \mathfrak p$ and $f_{\mathfrak p} = [k_{\mathfrak p} : \FF_\ell]$.  Let $S_{\mathfrak p}$ denote the set of embeddings $k_{\mathfrak p} \hookrightarrow \bar \FF_\ell$.  Define
\[
V_{\mathfrak p} = \bigtensor_{\tau \in S_{\mathfrak p}} \left ( \mathrm{det}^{a_{\tau}} \tensor_{k_{\mathfrak p}} \mathrm{Sym}^{b_{\tau} - 1} k_{\mathfrak p}^2 \right ) \tensor_{\tau} \bar \FF_\ell.
\]
Then the irreducible $\bar \FF_\ell$-representations of $G$ are of the form
\[
V = \bigtensor_{\mathfrak p | \ell} V_{\mathfrak p}. %
\]
Each factor $V_{\mathfrak p}$ of $V$ acts on the corresponding factor $GL_2(\mathcal O/ \mathfrak p)$ of $GL_2(\mathcal O / \ell \mathcal O)$.  

For computational purposes, we will think of $\mathrm{Sym}^{b_{\tau} - 1}(k_{\mathfrak p}^2)$ as the space of homogeneous polynomials of degree $b_{\tau} - 1$ in two variables with coefficients in $k_{\mathfrak p}$.  We define the left action of $\mathrm{GL}_2(\mathcal O)$ on this space as follows: For $g \in \mathrm{GL}_2(\mathcal O)$, we reduce $g$ modulo $\frak p$ and then the action is given by 
\[
	\bar g \cdot P(X, Y) = P(dX - bY, -cX + aY), 
\]
for $\bar g = \twodmatrix a b c d \in \mathrm{GL}_2(\mathcal O/\frak p)$.  This action is consistent with that used in Stein \cite[p.123]{stein}.  Note also that this definition of the action extends to any matrix $g \in M_2(\mathcal O)$, where $M_2(\mathcal O)$ denotes the semigroup of $2 \times 2$ matrices with entries in $\mathcal O$ and non-zero determinant.  This extension of the definition will be used in the definition of Hecke operators in Section \ref{hecke_ops}.

Given a Galois representation $\rho: G_K \rightarrow GL_2(\bar \FF_\ell)$ over a totally real field $K$, Buzzard, Diamond and Jarvis describe in detail a recipe for the predicted Serre weights of corresponding modular forms.  For each prime $\mathfrak p$ of $K$ such that $\mathfrak p | \ell$, they define a set $W_{\mathfrak p}(\rho)$ of $GL_2(k_{\mathfrak p})$-representations; the conjectural weight set $W(\rho)$ then consists of Serre weights of the form $V = \tensor_{\bar \FF_\ell} V_{\mathfrak p}$ where each $V_{\mathfrak p} \in W_{\mathfrak p}(\rho)$.  The sets $W_{\mathfrak p}(\rho)$ depend only on the local behaviour of $\rho$ at $\ell$ and are defined in different ways depending on whether $\rho$ restricted to the decomposition group $G_{k_{\mathfrak p}}$ is reducible or irreducible. See \cite[p. 17-27]{bdj} for a complete description of the weight recipe.

\subsection {Cohomological mod $\ell$ Forms over $K$} \label{defn_modp_forms}

Let $\Gamma$ be a congruence subgroup of $\mathrm{GL}_2(\mathcal O)$ of level $\frak n$.   We will be working with the following two congruence subgroups in particular:
\[
	\Gamma_0(\frak n) = \left \{ \twodmatrix a b c d \in \mathrm{GL}_2(\mathcal O) \;\mid\; \twodmatrix a b c d \equiv \twodmatrix \ast \ast 0 \ast \mod \; \frak n \right \}
\]
and
\[
	\Gamma_1(\frak n) = \left \{ \twodmatrix a b c d \in \mathrm{GL}_2(\mathcal O) \;\mid\; \twodmatrix a b c d \equiv \twodmatrix \ast \ast 0 1 \mod \; \frak n \right \}.
\]

Our definition of modular forms involves Hecke operators $T_{\frak q}$, which we will define in Section \ref{hecke_ops}.  

\begin{dfn}
	We define a \emph{cohomological mod $\ell$ form of level $\frak n$ and Serre weight $V$} to be a non-trivial cohomology class $f \in H^2(\Gamma, V)$ which is a simultaneous eigenvector for the Hecke operators $T_{\frak q}$ for all primes $\frak q$ such that $\frak q \nmid \ell \frak n$.    
\end{dfn}

\begin{dfn}
	Let $\rho: G_K \rightarrow \mathrm{GL}_2(\bar \FF_{\ell})$ be a continuous, irreducible representation and $V$ a Serre weight.  We say that $\rho$ is \emph{modular of weight $V$} if there is a non-zero modular form $f \in H^2(\Gamma, V)$ such that the eigenvalue $a_{\frak q}$ of $T_{\frak q}(f)$ is equal to $\mathrm{tr}(\rho(\mathrm{Frob}_{\frak q}))$ in $\bar \FF_{\ell}$ for all primes $\frak q \nmid \ell \frak n$.   
\end{dfn}

\section {Question} \label{question}

The main question of concern in this paper is the following:

\begin{quest}
	Let $K$ be an imaginary quadratic field and suppose 
	\[
		\rho: G_K \rightarrow \mathrm{GL}_2(\bar \FF_{\ell})
	\]
	is a continuous, irreducible representation.  Is it true that $\rho$ is modular of weight $V$ for every Serre weight $V$ in the BDJ conjectural weight set $W(\rho)$?
\end{quest}

The goal of this paper is to provide computational evidence in support of a positive answer to this question.  In Section \ref{modforms} below, we provide a modular symbols description of the cohomology space of modular forms in which we are interested.  Using this description the author wrote code in C, using the PARI library \cite{PARI2}, to compute modular forms.  In Section \ref{galreps}, we compute Galois representations over $K = \QQ(i)$ and then in Section \ref{tables} we provide tables of modular forms (computed with the author's code) corresponding to those Galois representations.  Further examples are given in \cite{torrey}.  

\begin{remark}
	In Serre's original conjecture, he requires that the representation $\rho$ be odd.  For a representation of $G_K$ where $K$ is an imaginary quadratic field, there is no odd/even distinction.  
\end{remark}

\section {Computing Modular Forms over $\QQ(i)$} \label{modforms}

One of the big challenges of producing computational evidence for Serre's conjecture is getting from the theoretical description of modular forms to a description that can be used in computations.   

In our setting, we want to compute the space $H^2(\Gamma, V)$ for some congruence subgroup $\Gamma$ and some Serre weight $V$.  Instead of computing this cohomology group directly, we use Borel-Serre duality to compute a homology group with coefficients in the Steinberg module.  The homology group is computationally friendly because we can express the Steinberg module in terms of modular symbols.

\subsection {Borel-Serre Duality} \label{borelserre_duality}

To apply Borel-Serre duality, we need a few assumptions.  We assume that $K$ is an imaginary quadratic field with class number 1 so that $\mathcal O$ is a PID and we assume that the order of the torsion in $\Gamma$ is invertible in $\overline{\FF}_{\ell}$.   In \cite{ash}, Ash describes the Steinberg module $St$ in terms of modular symbols for arbitrary dimension $n$. Here, we restrict to the case $n = 2$.  

\begin{dfn}\label{def_steinberg}
Let $R = \overline{\FF}_\ell$.  Consider the set of formal $R$-linear sums of symbols $[v] = [v_1, v_2]$ where the $v_i$ are unimodular columns in ${\mathcal O}^2$, i.e., $v_i = \twodcol a b $ with $\gcd(a, b) = 1$.  Mod out by the $R$-module generated by the following elements:
\begin{enumerate}
    \item $[v_2, v_1] + [v_1, v_2]$;
    \item $[v] = [v_1, v_2]$ whenever $\det(v) = 0$; and
    \item $[v_1, v_3] - [v_1, v_2] - [v_2, v_3]$,
\end{enumerate}
where the $v_i$ again run over all unimodular columns in $\mathcal O^2$. This quotient module is the \emph{Steinberg module $St$}.    We call the symbols $[v] = [v_1, v_2]$ \emph{modular symbols}.  
\end{dfn}

With the above assumptions, Borel-Serre duality \cite[p.482-483]{borelserre} gives an isomorphism
\[
    H^2(\Gamma, V) \overset{\sim}{\longrightarrow} H_0 (\Gamma, St \tensor V).
\]

Computing $H_0 (\Gamma, St \tensor V)$ amounts to computing the quotient of the Steinberg module $St$ by the relations $[v] - \gamma[v]$ for all $\gamma \in \Gamma$.

\subsection {An Algebraic Proposition} \label{alg_prop}

The following proposition is analogous to Proposition 4.3 in the doctoral thesis of Martin \cite[p.69]{martin}.  It will be used in Section \ref{manin_symbols} below to relate modular symbols to something called Manin symbols (which is what we will actually be computing).  First, we will need some notation.  We keep $R = \overline{\FF}_\ell$ as that is all we need here, though Proposition \ref{prop_algprop} will work for a general ring $R$.  In the following, we restrict to $K = \QQ(i)$.  To work with another imaginary quadratic field, one needs to prove some similar proposition that will depend on that particular field.  Define the following matrices in $\mathrm{GL}_2(\mathcal O)$:
\[
	J = \twodmatrix i 0 0 1, \qquad S = \twodmatrix 0 i 1 0, \qquad T = \twodmatrix 1 1 0 1, \qquad T' = \twodmatrix 1 0 1 1.
\]

\begin{prop}\label{prop_algprop}
	Consider the following homomorphism of left $R[\mathrm{GL}_2(\mathcal O)]$-modules:
	\[
	\begin{aligned}
		\Psi : R[\mathrm{GL}_2(\mathcal O)] &\longrightarrow R[\PP^1(K)] \\
		\sum_M u_M [M] &\longmapsto \sum_M u_M ([M(\infty)] - [M(0)]).
	\end{aligned}
	\]
	The kernel of $\Psi$ is equal to the left $R[\mathrm{GL}_2(\mathcal O)]$ ideal 
	\[
		\mathcal J = \langle [I] - [T] - [T'],\; [I] + [S],\; [I] - [J] \rangle.
	\]
\end{prop}
\begin{pf}
	First, we will show $\mathcal J \subseteq \ker(\Psi)$.  For this we simply evaluate $\Psi$ on $[I] - [T] - [T']$, on $[I] + [S]$ and on $[I] - [J]$.  We have
	\[
	\begin{aligned}
		\Psi([I] - [T] - [T']) &= [I(\infty)] - [I(0)] - [T(\infty)] + [T(0)] - [T'(\infty)] + [T'(0)] \\
			&= [\infty] - [0] - [\infty] + [1] - [1] + [0] \\
			&= 0
	\end{aligned}
	\]
	and
	\[
	\begin{aligned}
		\Psi([I] + [S]) &= [I(\infty)] - [I(0)] + [S(\infty)] - [S(0)] \\
			&= [\infty] - [0] + [0] - [\infty] \\ 
			&= 0
	\end{aligned}
	\]
	and
	\[
	\begin{aligned}
		\Psi([I] - [J]) &= [I(\infty)] - [I(0)] - [J(\infty)] + [J(0)] \\
			&= [\infty] - [0] - [\infty] + [0] \\ 
			&= 0.
	\end{aligned}
	\]
	
	The other direction, proving that $\ker(\Psi) \subseteq \mathcal J$, requires more work.  Let $W = \sum_M u_M [M]$ be a non-zero element of $\ker(\Psi)$.  Let $\mathcal L (W) \subseteq \PP^1(K)$ be the union of supports of $\sum_M u_M[M(\infty)]$ and $\sum_M u_M[M(0)].$  Furthermore, we define
	\[
		L(W) = \max_{\frac{\alpha}{\beta} \in \mathcal L(W)} (|\alpha|^2 + |\beta|^2)
	\]
	and
	\[
		m(W) = \left | \left \{ \frac{\alpha}{\beta} \in \mathcal L(W) : |\alpha|^2 + |\beta|^2 = L(W) \right \} \right |
	\]
	where we assume $(\alpha, \beta) = 1$.  We will use elements of the ideal $\mathcal J$ to write down a $W'$ congruent to $W$ modulo $\mathcal J$, but such that $L(W') \leq L(W)$.  Futhermore, if $L(W') = L(W)$ then we will have $m(W') < m(W)$.  Iterating this process, we will see that $W \in \mathcal J$.  
	
	Let $\alpha/\beta \in \mathcal L(W)$ be such that $|\alpha|^2 + |\beta|^2 = L(W)$.  Let $\delta$, $\gamma$ be elements of $\mathcal O$ such that $\alpha \gamma - \beta \delta = 1$ and such that $|\gamma| \leq |\beta|$ and $|\delta| \leq |\alpha|$.  Then the matrices of $\mathrm{GL}_2(\mathcal O)$ satisfying $M(\infty) = \alpha/\beta$ are of the form
	\[
		M = \twodmatrix{i^m \alpha}{i^n (\delta + k \alpha)}{i^m \beta}{i^n(\gamma + k \beta)}
	\]
	for $k \in \mathcal O$ and $m, n \in \{0, 1, 2, 3 \}$.  As $L(W) = |\alpha|^2 + |\beta|^2$, we see that for such a matrix $M$ to be in the support of $W$, we must have
	\[
		|\delta + k \alpha|^2 + |\gamma + k \beta|^2 \leq |\alpha|^2 + |\beta|^2.
	\]
	We will show first that for this to be true, we must have $|k| < 2$.  First, suppose $|k| \geq 2$.  Then 
	\[
	\begin{aligned}
		|\delta + k \alpha|^2 + |\gamma + k \beta|^2 &\geq (|k \alpha| - |\delta|)^2 + (|k \beta| - |\gamma|)^2 \\
			&\geq (2|\alpha| - |\delta|)^2 + (2|\beta| - |\gamma|)^2 \\
			&= (|\alpha| + (|\alpha| - |\delta|))^2 + (|\beta| + (|\beta| - |\gamma|))^2 \\
			&\geq |\alpha|^2 + |\beta|^2. \\
	\end{aligned}
	\]
	Note, furthermore, that equality is only possible if $|\alpha| = |\delta|$ and $|\beta| = |\gamma|$, but this contradicts the fact that $\alpha \gamma - \beta \delta = 1$.  Thus we must have $|k| < 2$.  Since $k \in \mathcal O$, the only possibilities are
	\begin{enumerate}
		\item $k = 0$;
		\item $k \in \mathcal O^{\ast}$; or
		\item $|k|^2 = 2$, i.e., $k \in \{ 1+i,\, 1-i,\, -1+i,\, -1-i \}$.
	\end{enumerate}
	
	We show why the assumption that $|\alpha| = |\delta|$ and $|\beta| = |\gamma|$ implies the determinant of $M$ cannot be in $\mathcal O^{\ast}$.  We have $\det(M) = i^{m+n} (\alpha \gamma - \beta \delta)$.  Under our assumption, we have $| \alpha \gamma | = |\beta \delta |$, so we need to show that we cannot have $x - y \in \mathcal O^{\ast}$ if $|x| = |y|$.  Suppose this is the case and let $x = a+bi$ and $y = c+di$, so that $x - y = (a-c) + (b-d)i$.  Then either $a = c$ and $b = d \pm 1$ or $b = d$ and $a = c \pm 1$.  Without loss of generality, assume $a = c$ and $b = d \pm 1$.  Then 
		\[
		\begin{aligned}
			 |x|^2 &= |y|^2 \\
			\Rightarrow \qquad a^2 + b^2 &= c^2 + d^2 \\
			\Rightarrow \qquad (d \pm 1)^2 &= d^2, \\
		\end{aligned}
		\]
		giving a contradiction.
	
	Using the fact that we are taking the quotient by the elements $I - J$ and $I + S$ of $\mathcal J$, we may assume the matrix $M$ has the form
	\[
		M_k = \twodmatrix{\alpha}{\delta + k \alpha}{\beta}{\gamma + k \beta}.
	\]
	
	Similarly, we may assume that the only matrices $N$ in the support of $W$ such that $N(0) = \alpha/\beta$ are those of the form
	\[
		N_j = \twodmatrix{j \alpha - \delta}{\alpha}{j \beta - \gamma}{\beta},
	\]
	where $|j|^2 \in \{0,\; 1,\; 2 \}$  and $|j \alpha - \delta|^2 + |j \beta - \gamma|^2 \leq |\alpha|^2 + |\beta|^2$.  Each of the $N_j$ can be replaced, using $S$ and $J$, with $-M_{-j}$ as follows:  Note that $I + JS \in \mathcal J$ since $I + JS = (I - J) + J(I + S)$.  Then 
	\[
		N_j - N_j(I + JS) = - N_j JS = - M_{-j}.
	\]
	
	Let $W'$ denote the element $W \in \ker (\Psi)$ with the above modifications.  So now all matrices $M$ in the support of $W'$ with $M(\infty) = \alpha/\beta$ are of the form
	\[
		M_k = \twodmatrix{\alpha}{\delta + k \alpha}{\beta}{\gamma + k \beta},
	\]
	with
	\[
		|\delta + k \alpha|^2 + |\gamma + k \beta|^2 \leq |\alpha|^2 + |\beta|^2.
	\]
	Furthermore, we have no matrices $N$ in the support of $W'$ with $N(0) = \alpha/\beta$.  
	
	The strategy from here is as follows:  We will first consider the case where $|k|^2 = 2$, and we will replace such matrices $M$ with two matrices $M'$ and $M''$, where $M'$ will have the same form as $M$ above but with $|k| = 1$ and $M''$ will be such that $M''(\infty) \neq \alpha/\beta \neq M''(0)$.  These matrices will also satisfy $m(M) = m(M' - M'')$, so that the number of $\alpha/\beta \in \mathcal L(W)$ will not increase.  The next step is to work with the case $|k| = 1$.  We will replace each $M$ of this sort again with two matrices $M'$ and $M''$.  One of these will be in the form of $M_k$ but with $k = 0$.  The other will again be such that $M''(\infty) \neq \alpha/\beta \neq M''(0)$.  As in the $|k|^2 = 2$ case, these matrices will satisfy $m(M) = m(M' - M'')$.  After these steps, the only matrix $M$ left in $W'$ with $M(\infty) = \alpha/\beta$ is $M_0$, and there are no matrices left in $W'$ with $M(0) = \alpha/\beta$.  Since $W' \in \ker (\Psi)$, we must have the coefficient $u_{M_0}$ of $M_0$ in $W'$ equal to $0$.  We thus decrease $m(W')$ by at least one.  
	
	To deal with the case $|k|^2 = 2$, we use the following lemma.
	\begin{lemma}\label{clm_algprop_absvalk2}
		Suppose $k \in \mathcal O$ is such that $|k|^2 = 2$, i.e., $k \in \{ 1+i,\, 1-i,\, -1+i,\, -1-i \}$ and write $k = k_1 + k_2i$.  If 
		\[
			|\delta + k \alpha|^2 + |\gamma + k \beta|^2 \leq |\alpha|^2 + |\beta|^2,
		\]
		then either $|\delta + k_1 \alpha|^2 + |\gamma + k_1 \beta|^2$ or $|\delta + k_2i \alpha|^2 + |\gamma + k_2i \beta|^2$ is strictly less than $|\alpha|^2 + |\beta|^2$.
	\end{lemma}
	\begin{pf}
	The proof is elementary (see \cite[p. 67-69]{torrey}).
	\end{pf}

	For the first step, i.e. $|k|^2 = 2$, we can apply Lemma \ref{clm_algprop_absvalk2} as follows: Let $t_1 \in \{ k_1, k_2i \}$ be such that $|\delta + t_1 \alpha|^2 + |\gamma + t_1 \beta| < |\alpha|^2 + |\beta|^2$ and let $t_2$ be the other element of $\{ k_1, k_2i \}$.  Then 
	\[
		\twodmatrix{t_2 \alpha}{\delta + t_1 \alpha}{t_2 \beta}{\gamma + t_1 \beta} T = \twodmatrix{t_2 \alpha}{\delta + t_1 \alpha}{t_2 \beta}{\gamma + t_1 \beta} \twodmatrix 1 1 0 1 = \twodmatrix{t_2 \alpha}{\delta + k \alpha}{t_2 \beta}{\gamma + k \beta} 
	\]
	and 
	\[
		\twodmatrix{t_2 \alpha}{\delta + t_1 \alpha}{t_2 \beta}{\gamma + t_1 \beta} T' = \twodmatrix{t_2 \alpha}{\delta + t_1 \alpha}{t_2 \beta}{\gamma + t_1 \beta} \twodmatrix 1 0 1 1 = \twodmatrix{\delta + k \alpha}{\delta + t_1 \alpha}{\gamma + k \beta}{\gamma + t_1 \beta}. 
	\]
	Since $t_2 = i^m$ for some integer $m$, we can apply $I - J$ to $M$ until we have 
	\[
		M = \twodmatrix{t_2 \alpha}{\delta + k \alpha}{t_2 \beta}{\gamma + k \beta}.
	\]
	Defining
	\[
		M' = \twodmatrix{t_2 \alpha}{\delta + t_1 \alpha}{t_2 \beta}{\gamma + t_1 \beta}
	\]
	and 
	\[
		M'' = \twodmatrix{\delta + k \alpha}{\delta + t_1 \alpha}{\gamma + k \beta}{\gamma + t_1 \beta},
	\]
	we then have 
	\[
	\begin{aligned}
		M' (I - T - T') &= \twodmatrix{t_2 \alpha}{\delta + t_1 \alpha}{t_2 \beta}{\gamma + t_1 \beta} ( I - T - T') \\
			&= \twodmatrix{t_2 \alpha}{\delta + t_1 \alpha}{t_2 \beta}{\gamma + t_1 \beta} - \twodmatrix{t_2 \alpha}{\delta + k \alpha}{t_2 \beta}{\gamma + k \beta} - \twodmatrix{\delta + k \alpha}{\delta + t_1 \alpha}{\gamma + k \beta}{\gamma + t_1 \beta} \\
			&= M' - M - M'',
	\end{aligned}
	\]
	so, modulo $\mathcal J$, we may replace $M$ with $M' - M''$.  Note that $M'$ can be replaced with 
	\[
		M_{k'} = \twodmatrix{\alpha}{\delta + k' \alpha}{\beta}{\gamma + k' \beta},
	\]
	with $|k'|^2 = 1$, a case we will treat shortly.  Also, we have $L(M'') \leq L(M)$ and $M''(\infty) \neq \alpha/\beta \neq M''(0)$.  Recall that when $L(W') = L(W)$, we need $m(W') < m(W)$.  In this step, we are replacing one matrix $M$ with two, $M'$ and $M''$, but since $|\delta + t_1 \alpha|^2 + |\gamma + t_1 \beta| < |\alpha|^2 + |\beta|^2$, we have $m(M) = m(M' - M'')$.  (Note that $M''(\infty) = (\delta + k \alpha)/(\gamma + k \beta)$ was already in $\mathcal L(W')$ as $M''(\infty) = M(0)$.)  We are not decreasing $m(W)$ in this step, but we are not increasing it either.  In later steps we will obtain a decrease in either $L(W')$ or $m(W')$.  
	
	So now the only matrices $M$ remaining in the support of $W$ with $M(\infty) = \alpha/\beta$ are those of the form
	\[
		M_k = \twodmatrix{\alpha}{\delta + k \alpha}{\beta}{\gamma + k \beta},
	\]
	where either $k = 0$ or $|k|^2 = 1$ and $|\delta + k \alpha|^2 + |\gamma + k \beta|^2 \leq |\alpha|^2 + |\beta|^2$ (by assumption for those starting out with $|k|^2 = 1$ and by design for those coming initially from matrices with $|k|^2 = 2$).  
	
	Now we will use elements of $\mathcal J$ to eliminate the $M_k$ for $|k| = 1$.  In particular, we will use $I - T - T'$.   We can first use repeated applications of $I - J$ to replace $M_k$ with $\tilde M_k$ with left hand column given by the transpose of $({-k \alpha} \;\; {-k \beta})$, i.e.,    
	\[
		\tilde M_k = \twodmatrix{-k\alpha}{\delta + k \alpha}{-k\beta}{\gamma + k \beta}.
	\]
	Defining
	\[
		M' = \twodmatrix{-k\alpha}{\delta}{-k\beta}{\gamma} \;\;\mathrm{ and }\;\; M'' = \twodmatrix{\delta}{\delta + k \alpha}{\gamma}{\gamma + k \beta},
	\]
	we then have 
	\[
	\begin{aligned}
		\tilde M_k (I - T - T') &= \twodmatrix{-k\alpha}{\delta + k \alpha}{-k\beta}{\gamma + k \beta} ( I - T - T') \\
			&= \twodmatrix{-k\alpha}{\delta + k \alpha}{-k\beta}{\gamma + k \beta} - \twodmatrix{-k\alpha}{\delta}{-k\beta}{\gamma} - \twodmatrix{\delta}{\delta + k \alpha}{\gamma}{\gamma + k \beta} \\
			&= \tilde M_k - M' - M''.
	\end{aligned}
	\]
	
	Thus we can replace $\tilde M_k$ with $M' + M''$, where $M'$ can be replaced with $M_0$ and $M''$ will satisfy:
	\begin{enumerate}
		\item $M''(\infty) \neq \alpha/\beta \neq M''(0)$ and
		\item $L(M'') \leq L(W)$.
	\end{enumerate}
	Furthermore we have $m(\tilde M_k) = m(M' + M'')$.  Note that $M''(0) = (\delta + k \alpha)/(\gamma + k \beta)$ was already in $\mathcal L(W')$ as $M''(0) = \tilde M_k(0)$.
	
	Finally the only matrix in the support of $W'$ with either $M(\infty)$ or $M(0)$ equal to $\alpha/\beta$ is $M_0$.  Thus we must have $u_{M_0} = 0$ and so $m(W')$ is decreased by at least one.  	
	\end{pf}

\subsection {Manin Symbols} \label{manin_symbols}

We will follow the approach of Wiese \cite[p.6-9]{wiese_hecke} in using Proposition \ref{prop_algprop} from Section \ref{alg_prop} to relate modular symbols to Manin symbols (which we will define at the end of this subsection).  Manin symbols provide us with an explicit, computationally friendly description of the homology we wish to compute.  Proposition \ref{prop_modtomanin} below gives us the first step in the transition from modular symbols to Manin symbols.  We define another matrix:
\[
	L = \twodmatrix {1}{-1}{1}{0}.
\]

In the following we will use the notation $\alpha/\beta$ for the unimodular column  $\twodcol{\alpha}{\beta}$.  In particular, we will use $0$ to denote $\twodcol 0 1$ and $\infty$ to denote $\twodcol 1 0 $.  

\begin{prop}\label{prop_modtomanin}
	The following homomorphism of $R$-modules is an isomorphism:
	\[
	\begin{aligned}
		\Phi: R[\mathrm{GL}_2(\mathcal O)]/\mathcal I &\longrightarrow St \\
		M &\longmapsto [ M (0), M (\infty) ]
	\end{aligned}
	\]
	where 
	\[
		\mathcal I = R[\mathrm{GL}_2(\mathcal O)](I - J) + R[\mathrm{GL}_2(\mathcal O)](I + S) + R[\mathrm{GL}_2(\mathcal O)](I + L + L^2).
	\]
\end{prop}
\begin{pf}
	To prove that $\Phi$ is surjective, first note that, in the Steinberg module $St$, we have 
	\[
	\begin{aligned}[] 
		[ v_1, v_2 ] &= [v_1, 0] + [0, v_2] \\
			&= -[0, v_1] + [0, v_2], \\
	\end{aligned}
	\]
	so it suffices to show $[0, \alpha/\beta]$ is in the image of $\Phi$, where $\alpha/\beta$ is any unimodular column in $\mathcal O_2$.  
	
	We use continued fractions to write a given modular symbol of the form $[ 0, \alpha/\beta ]$ as a finite sum of symbols of the form $[\gamma(0), \gamma(\infty) ]$ with $\gamma \in \mathrm{SL}_2(\mathcal O)$.  For this algorithm, we rely on the fact that the field $K$ is Euclidean.  We set $r_0 = \alpha/\beta$ and $r_n = 1/(r_{n-1} - a_{n-1})$ and define $a_n = \mathrm{floor}(r_n)$, where we define the floor function on $\QQ(i)$ as follows: $\mathrm{floor}(r + s i)$ for $r, \; s \in \QQ$ is equal to $a + bi$ with $a$ and $b$ the least integers such that $|a - r| \leq 1/2$ and $|b - s| \leq 1/2$.
	
	We then define the convergents of the continued fractions as follows:
	\[
	\begin{aligned}
		p_{-2} = 0 &\qquad q_{-2} = 1 \\
		p_{-1} = 1 &\qquad q_{-1} = 0 \\
		p_{n} = a_n p_{n-1} + p_{n-2} &\qquad q_{n} = a_n q_{n-1} + q_{n-2} 
	\end{aligned}
	\]
	so that 
	\[
		\frac{p_n}{q_n} = a_0 + \frac{1}{a_1 + \frac{1}{a_2 + \frac{1}{\cdots + \frac{1}{a_n}}}}
	\]
	and for some $k$ we have
	\[
		\frac{\alpha}{\beta} = \frac{p_k}{q_k}.
	\]
	As in continued fractions for $\ZZ$, we have
	\[
		p_n q_{n-1} - p_{n-1} q_n = (-1)^{n+1}.
	\]
	We can now rewrite the modular symbol $[ 0, \alpha/\beta ]$ as
	\[
		\left [ 0, \; \frac{\alpha}{\beta} \right ] = \sum_{n=-1}^k \left [ \frac{p_{n-1}}{q_{n-1}}, \frac{p_n}{q_n} \right ] = \sum_{n=-1}^k [ \gamma_n(0), \gamma_n(\infty) ] = \sum_{n=-1}^k \Phi(\gamma_n),
	\]
	where
	\[
		\gamma_n = \twodmatrix{(-1)^{n+1}p_n}{p_{n-1}}{(-1)^{n+1}q_n}{q_{n-1}}.
	\]
	
	We now show that the kernel of $\Phi$ is equal to 
	\[
		\mathcal I = R[\mathrm{GL}_2(\mathcal O)](I - J) + R[\mathrm{GL}_2(\mathcal O)](I + S) + R[\mathrm{GL}_2(\mathcal O)](I + L + L^2).
	\]
	We start by showing that $\ker(\Phi) = \ker(\Psi)$ where $\Psi$ is the map in Proposition \ref{prop_algprop}.  Define another map $\pi$ by
	\[
	\begin{aligned}
		\pi: St &\longrightarrow R[\PP^1(K)] \\
		[v_1, v_2] &\longmapsto v_2 - v_1.
	\end{aligned}
	\]
	Then we have $\Psi = \pi \circ \Phi$, so certainly $\ker(\Phi) \subseteq \ker(\Psi)$.  To show the other inclusion, suppose $\sum_{M} u_M M \in \ker (\Psi)$, i.e., 
	\[
		\Psi \left ( \sum_M u_M M \right ) = \sum_M u_M M(0) - \sum_M u_M M(\infty) = 0.
	\]
	Applying $\Phi$ instead of $\Psi$ to this element of $\ker(\Psi)$ and using the first of the relations in Definition \ref{def_steinberg} (definition of $St$), we then have
	\[
	\begin{aligned}
		\Phi \left ( \sum_M u_M M \right ) &= \sum_M u_M \left [ M (0), M (\infty) \right ] \\
			&= \sum_M u_M \left [ M (0), \infty \right ]
				+ \sum_M u_M \left [ \infty, M (\infty) \right ] \\
			&= \sum_M u_M \left [ M (0), \infty \right ]
				- \sum_M u_M \left [ M (\infty), \infty \right ] \\
			&= 0,
	\end{aligned}
	\]
	and so $\ker(\Psi) \subseteq \ker(\Phi)$.
	
	Finally, we need to show that $\ker(\Phi)$ can be written in the form claimed, i.e., that
	\[
		\ker(\Phi) = R[\mathrm{GL}_2(\mathcal O)](I - J) + R[\mathrm{GL}_2(\mathcal O)](I + S) + R[\mathrm{GL}_2(\mathcal O)](I + L + L^2).
	\]
	Currently, we have that 
	\[
		\ker(\Phi) = R[\mathrm{GL}_2(\mathcal O)](I - J) + R[\mathrm{GL}_2(\mathcal O)](I + S) + R[\mathrm{GL}_2(\mathcal O)](I - T - T').
	\]
	First, note that in
	\[
		R[\mathrm{GL}_2(\mathcal O)](I - J) + R[\mathrm{GL}_2(\mathcal O)](I + S)
	\]
	we have both $I + \tilde S$ and $I - \tilde J$ where 
	\[
		\tilde S  = JS = \twodmatrix{0}{-1}{1}{0} \quad \mathrm{and} \quad \tilde J = J^2 = \twodmatrix{-1}{0}{0}{1}
	\]
	as 
	\[
		I + \tilde S = (I - J) + J(I + S) 
	\]
	and 
	\[
		I - \tilde J = (I - J) + J(I - J). 
	\]
	We then also have $I + \tilde J \tilde S \tilde J$ as
	\[
		I + \tilde J \tilde S \tilde J = (I - \tilde J) + \tilde J (I + \tilde S) - \tilde J \tilde S (I - \tilde J).
	\]
	
	Now, to see that the two forms of the kernel are the same, we write
	\[
	\begin{aligned}
		(I - T - T') + T(I + \tilde S) + T'(I + \tilde S) &= I + T \tilde S + T' \tilde S \\
		&= I + L + L^2
	\end{aligned}
	\]
	and 
	\[
	\begin{aligned}
		(I + L + L^2) - (L + L^2)(I + \tilde J \tilde S \tilde J) &= I - L\tilde J \tilde S \tilde J - L^2 \tilde J \tilde S \tilde J \\
		&= I - T - T'.
	\end{aligned}
	\]
\end{pf}

Our discussion thus far is only sufficient for computing weight two modular forms.  We now extend this so that we can compute higher weight forms and we also incorporate the $\Gamma$ relations $[v] - \gamma [v]$.  

Let $R = \bar \FF_{\ell}$ and let $V$ be a Serre weight with left $R[\mathrm{GL}_2(\mathcal O)]$ action, as defined in Section \ref{serreweights}.   We consider the module of left $\Gamma$-coinvariants $(R[\mathrm{GL}_2(\mathcal O)] \tensor_R V)_{\Gamma}$, where $\Gamma$ acts diagonally on the left and we have the natural right $R[\mathrm{GL}_2(\mathcal O)]$ action, i.e., $(h \tensor v)g = (hg \tensor v)$.  

In the following theorem and subsequent proposition we use Proposition \ref{prop_modtomanin} to write the homology group we wish to compute in terms of Manin symbols.

\begin{theorem}
	Let $N$ denote the $R$-module $(R[\mathrm{GL}_2(\mathcal O)] \tensor_R V)_{\Gamma}$ as described above.  Then the following sequence of $R$-modules is exact:
	\[
		0 \rightarrow N(I - J) + N(I + S) + N(I + L + L^2) \rightarrow N \rightarrow H_0(\Gamma, St \tensor V) \rightarrow 0.
	\]
\end{theorem}
\begin{pf}
	Proposition \ref{prop_modtomanin} gives the exact sequence
	\[
		0 \rightarrow \mathcal I \rightarrow R[\mathrm{GL}_2(\mathcal O)] \rightarrow St \rightarrow 0
	\]
	where
	\[
		\mathcal I = R[\mathrm{GL}_2(\mathcal O)](I - J) + R[\mathrm{GL}_2(\mathcal O)](I + S) + R[\mathrm{GL}_2(\mathcal O)](I + L + L^2).
	\]
	Let $N' = R[\mathrm{GL}_2(\mathcal O)] \tensor_R V$.  Since $V$ is a free $R$-module,  the following sequence of $R[\Gamma]$-modules is also exact:
	\[
		0 \rightarrow N'(I - J) + N'(I + S) + N'(I + L + L^2) \rightarrow N' \rightarrow St \tensor V \rightarrow 0.
	\]
	We then need only take $\Gamma$-coinvariants to achieve the desired exact sequence.  
\end{pf}

We need one more step to get to the module we will actually be computing.  This is the content of the following proposition.
\begin{prop}\label{prop_manin_mod}
	Let $X$ denote the $R$-module $R[\Gamma \backslash \mathrm{GL}_2(\mathcal O)] \tensor_R V$ with right $\mathrm{GL}_2(\mathcal O)$-action given by $(\Gamma h \tensor v) g = \Gamma h g \tensor g^{-1} v$.  Then there is a right $R[\mathrm{GL}_2(\mathcal O)]$-module isomorphism between $N = (R[\mathrm{GL}_2(\mathcal O)] \tensor_R V)_{\Gamma}$ and $X$. 
\end{prop}
\begin{pf}
	The isomorphism is simply given by $g \tensor v \mapsto g \tensor g^{-1}v$.  
\end{pf}

The $R$-module $X = R[\Gamma \backslash \mathrm{GL}_2(\mathcal O)] \tensor_R V$ is the module of Manin symbols and is the basic module we will use in computations.  For $\Gamma = \Gamma_0(\frak n)$, the coset representatives of $\Gamma \backslash \mathrm{GL}_2(\mathcal O)$ are in one-to-one correspondence with $\PP^1(\frak n)$, the projective line over $\mathcal O/\frak n$.  We use the notation $(c:d)$ with $c,\; d \in \mathcal O$ to denote such a coset representative.

\subsection {$\Gamma_1(\frak n)$ and Characters} \label{gamma_1_chars}

The treatment in Section \ref{manin_symbols} is sufficient for dealing with $\Gamma = \Gamma_0(\frak n)$.  We also want to compute modular forms for $\Gamma_1(\frak n)$.  We will follow the approach used by Wiese \cite{wiese_hecke} for the $\Gamma_1(\frak n)$ case.

In the $\Gamma_1(\frak n)$ case, we will use a character
\[
	\varepsilon : \Gamma_0(\frak n) \onto \Gamma_1(\frak n) \backslash \Gamma_0(\frak n) \overset{\sim}{\rightarrow} (\mathcal O/\frak n)^{\ast} \rightarrow \bar \FF_{\ell}^{\ast}.
\]
We define a slight variation on the weight module $V$, which takes into account the action of the character $\varepsilon$.  This we define as 
\[
	V^{\varepsilon} = V \tensor_{\bar \FF_{\ell}} \bar \FF_{\ell}^{\varepsilon},
\]
where $\bar \FF_{\ell}^{\varepsilon}$ denotes a copy of $\bar \FF_{\ell}$ with action of $\Gamma_0(\frak n)$ by $\varepsilon^{-1}$.  

Computing cohomological mod $\ell$ modular forms for $\Gamma_1(\frak n)$ with character $\varepsilon$ then amounts to computing simultaneous eigenvectors for the Hecke operators on 
\[
	N^{\varepsilon} = \left ( \left ( R[\mathrm{GL}_2(\mathcal O)] \tensor V^{\varepsilon} \right ) _{\Gamma_1(\frak n)} \right )_{\Gamma_1(\frak n) \backslash \Gamma_0(\frak n)} ,
\]
modulo the relations used in Proposition \ref{prop_modtomanin}.  Here $\Gamma_1(\frak n) \backslash \Gamma_0(\frak n)$ is acting diagonally on the left on $\left ( R[\mathrm{GL}_2(\mathcal O)] \tensor V^{\varepsilon} \right )_{\Gamma_1(\frak n)} $. 
\begin{prop}\label{prop_gamma1_msymbols}
	Consider the $R$-module
	\[
		X = R[\Gamma_0(\frak n) \backslash \mathrm{GL}_2(\mathcal O)] \tensor V \tensor {\bar \FF_{\ell}}^{\varepsilon},
	\]
	where $\mathrm{GL}_2(\mathcal O)$ acts on the right by $(h \tensor v \tensor r)g = (hg \tensor g^{-1} v \tensor r)$ and $\Gamma_1(\frak n) \backslash \Gamma_0(\frak n)$ acts on the left by $g(h \tensor v \tensor r) = (gh \tensor v \tensor \varepsilon(g) r)$.   We have
	\[
		N^{\varepsilon} \cong X.
	\]
\end{prop}
\begin{pf}
	First we apply the isomorphism 
	\[
		 \left ( \left ( R[\mathrm{GL}_2(\mathcal O)] \tensor V^{\varepsilon} \right )_{\Gamma_1(\frak n)}  \right )_{\Gamma_1(\frak n) \backslash \Gamma_0(\frak n)} \cong   \left ( R[\mathrm{GL}_2(\mathcal O)] \tensor V^{\varepsilon} \right )_{\Gamma_0(\frak n)} . 
	\]
	Then we apply the isomorphism of Proposition \ref{prop_manin_mod}.  
\end{pf}

\subsection {Hecke Operators} \label{hecke_ops}

In this section, we define Hecke operators on the space of modular symbols $H_0(\Gamma, St \tensor V \tensor {\bar \FF_{\ell}}^{\varepsilon})$, with $\Gamma = \Gamma_0(\frak n)$.  To compute eigenvalues, we convert Manin symbols to modular symbols, compute the action of the Hecke operators there, and then convert back to Manin symbols.  We use the results of Section \ref{manin_symbols} to convert back and forth.  

Let $\frak p$ be a prime ideal of $\mathcal O$ that is relatively prime to the level $\frak n$ and let $\pi$ be a generator for $\frak p$.  We define a set $\Delta_{\frak p} \subset \mathrm{GL}_2(K)$ by
\[
	\Delta_{\frak p} = \left \{ \twodmatrix a b c d \in M_2(\mathcal O) \; : \; (ad - bc)\mathcal O = \frak p, \; \twodmatrix a b c d \equiv \twodmatrix{u}{\ast}{0}{\pi} \mod \; \frak n \right \},
\]
where $u \in \mathcal O^{\ast}$.  Let $\sigma_a \in \mathrm{SL}_2(\mathcal O)$ be defined by $\sigma_a \equiv \twodmatrix{1/a}{0}{0}{a} \mod \; \frak n$ for $a \in \{1, \; \pi \}$. Using the fact that $K$ is Euclidean, one can easily show that
\[
	\Gamma_1(\frak n) \Delta_{\frak p} = \Delta_{\frak p} \Gamma_1(\frak n) = \Delta_{\frak p}
\]
and, taking $a \in \{1, \; \pi \}$ and letting $x$ run over representatives of $\mathcal O/(\pi/a)$, $\Delta_{\frak p}$ can be written as a disjoint union as
\[
	\Delta_{\frak p} = \bigcup_{a,x} \Gamma_1(\frak n) \cdot \sigma_a \twodmatrix{a}{x}{0}{\pi/a}.
\]
Furthermore, since we take coinvariants via the character action on $\Gamma_1(\frak n) \backslash \Gamma_0(\frak n)$ and $\sigma_a \in \Gamma_0(\frak n)$, we may define the Hecke operator $T_{\frak p}$ by 
\[
\begin{aligned}
	T_{\frak p} \left ( [ v_1, v_2 ] \tensor v \tensor r \right ) = &\varepsilon(\pi) \twodmatrix {\pi}{0}{0}{1} \left ( [ v_1, v_2 ] \tensor v \tensor r \right ) \\
	&+ \sum_{x \; \mod \; \pi} \twodmatrix {1}{x}{0}{\pi} \left ( [ v_1, v_2 ] \tensor v \tensor r \right ). 
\end{aligned}
\]
The Hecke operator $T_{\frak p}$ is well-defined since $\Gamma_1(\frak n) \Delta_{\frak p} = \Delta_{\frak p} \Gamma_1(\frak n)$.   It is easy to see that the Hecke operator $T_{\frak p}$ is independent of the choice of generator $\pi$ using the relations in Proposition \ref{prop_algprop}.

This Hecke action on homology is compatible with the usual Hecke action on cohomology (mentioned in Section \ref{defn_modp_forms}) as shown in \cite[p.407]{agm}.

\section {Examples of Galois Representations} \label{galreps}

In this section, we compute examples of Galois representations.  These examples come from two sources: polynomials and class field theory.    Further examples from these sources and examples from elliptic curves can be found in \cite{torrey}.

Much of the data in this section was computed with the mathematical software systems Magma \cite{magma}, PARI/GP \cite{PARI2} and Sage \cite{sage}.

\subsection {Torsion} \label{torsion}

Early in section \ref{borelserre_duality} we assumed that the torsion in $\Gamma$ is invertible in the commutative ring $R$.  In our case, $\Gamma$ is a congruence subgroup of $GL_2(\mathcal O)$ and $R = \bar \FF_{\ell}$.  This assumption about torsion is not a strong assumption since we can make sure that, as long as the level $\frak n$ is large enough, the congruence subgroup $\Gamma$ will be torsion free.  In our examples, we have $K = \QQ(i)$ and $\Gamma = \Gamma_1(\frak n)$ for some level $\frak n$.  Now suppose $A \in GL_2(K) \setminus K^\ast$ is a torsion element of prime power, i.e., there is some prime $p$ such that $A^p = I$.  Now consider $K[A]$.  We have
\[
	K[A] \cong K[x]/(x^2 - Tx + D),
\]
where $T$ is the trace of $A$ and $D$ is the determinant of $A$.  Since $A \not \in K^\ast$, we may assume $x^2 - Tx + D$ is irreducible and so $K[A]$ is a field and a quadratic extension of $K$.  Since $A$ is a $p^{\mathrm{th}}$ root of unity, we also have $K[A] \cong K[\zeta_p]$.  Then we have a quadratic extension of $K$ of the form $K[\zeta_p]$.  Since $K = \QQ(i)$, this implies that $p$ is $2$ or $3$.  We will not use $\ell = 2$ for any of our examples (since $2$ is ramified in $K$ and hence not covered by the BDJ conjecture), but we do have some examples with $\ell = 3$, so consider $p = 3$.  Then in the polynomial above we have $T = -1$ and $D = 1$.  The matrix $A$ is in the congruence subgroup $\Gamma_1(\frak n)$, so we have 
\[
	A = \twodmatrix a b c d \equiv \twodmatrix {\ast}{\ast}{0}{1} \mod \; \frak n.
\]
Since the determinant $D = 1$, this implies that $a \equiv 1 \; \mod \; \frak n$.  Then $T = -1$ implies $2 \equiv -1 \; \mod \; \frak n$, which in turn implies that $\frak n \mid 3 \mathcal O$.  None of our examples have $\frak n \mid 3 \mathcal O$, and so all of our examples satisfy the condition on torsion in $\Gamma$.

\subsection {Representation from a Polynomial} \label{exsgalreps_polys}

In this section we will examine an example of a Galois representation arising from a polynomial.  We will determine the level, character, coefficients and predicted weights for this representation.  

In \cite[p 117]{figueiredo}, Figueiredo gives three examples of $A_4$ representations.  The first of these examples comes from the polynomial
\[
	x^4 - 7x^2 - 3x + 1 \qquad \mathrm{disc } = 3^2 \times 61^2.
\]
Figueiredo considers the mod $3$ representation arising from this polynomial using the isomorphism $A_4 \cong \mathrm{PSL}_2(\FF_3)$.  He shows that there must be a lift of this representation to $\mathrm{GL}_2(\bar \FF_3)$.  I will instead compute representations directly from the $\hat{A}_4$ extension, where $\hat{A}_4$ is a double cover of $A_4$, isomorphic to $\mathrm{SL}_2(\FF_3)$.  From the database of Kl\"{u}ners and Malle \cite{db_km}, we find the polynomial giving the $\hat{A}_4$ extension of Figueiredo's polynomial:
\[
	x^8+3x^7-11x^6-9x^5+21x^4+9x^3-11x^2-3x+1
\]
and let $L$ denote the Galois closure of this polynomial.

We take $\ell = 3$.  There is only one irreducible $2$-dimensional mod $3$ representation
\[
	\rho : G_K \rightarrow \mathrm{GL}_2(\bar \FF_3)
\]
factoring through $\mathrm{Gal}(LK/K)$.  We get this representation by taking the base change to $K = \QQ(i)$ from the representation $\rho_{\QQ}$ of $G_{\QQ}$, which we get simply by restricting $G_{\QQ}$ to $\mathrm{Gal}(L/\QQ)$ and then applying the following isomorphisms and inclusion:
\[
	\mathrm{Gal}(L/\QQ) \cong \hat{A}_4 \cong \mathrm{SL}_2(\FF_3) \hookrightarrow \mathrm{GL}_2(\FF_3).
\]

Note that this representation $\rho_{\QQ}$ is \emph{even}.  If we start with a representation $\rho_{\QQ}$ over $\QQ$ which is \emph{odd}, i.e. $\det \rho(\sigma_{\infty}) = -1$, then we already know $\rho_{\QQ}$ is modular by Serre's conjecture.  Instead, we start with an even representation $\rho_{\QQ}$, and look at the base change to $K$.  

To compute the level of $\rho_K$, we use the description and notation in \cite{serre2}.   We look at all primes dividing the discriminant of the polynomial, except $\ell$.  For each such prime $p$, we compute $n(p, \rho|_{\QQ})$.  In the examples we compute we will have $p$ unramified in $K$ for all primes $p$ dividing the level $N$ of the representation $\rho|_{\QQ}$.  When this is the case, we have $n(\frak p, \rho|_K) = n(p, \rho|_{\QQ})$ for primes $\frak p$ lying above $p$.   The corresponding level to check will be 
\[	
	\frak n = \prod \frak p^{n(\frak p, \rho)}.
\]
We know that for any prime $p$ which is tamely ramified in $L$, we have $n(p, \rho) = \mathrm{dim}(V/V_0)$, where $V_0$ is the subspace of $V$ fixed by the inertia group $I_{\frak P}$ for some prime $\frak P$ of $L$ lying above $p$.  If $p$ is wildly ramified in $L$, the determination of $n(p, \rho)$ requires further analysis, but this does not occur in any of our examples.  

For the level in this example, we need only compute $n(61, \rho)$.  We have the ramification index $e_{61} = 3$ in $L$, and so $61$ is tamely ramified and $n(61, \rho) = \dim (V/V_0)$.  We check that $\rho$ restricted to the subgroup of order $3$ fixes a $1$-dimensional subspace of $V$, so $n(61, \rho) = 1$ and the level of $\rho$ is $\frak n = 61$.  

Since the image of $\rho$ is $\mathrm{SL}_2(\FF_3)$, we see that $\det(\rho)$ is trivial, and so the character of $\rho$ is trivial as well.  

We want to compute the sequence $\{ a_{\mathfrak q} \}$ associated to the representation $\rho$, i.e. we want to compute $\mathrm{tr}(\rho(\mathrm{Frob}_{\mathfrak q}))$ for each prime $\mathfrak q$ of $\mathcal O$ where $\frak q \nmid \ell \frak n$.  Let $\mathfrak Q$ be a prime of $LK$ lying over $\mathfrak q$.  Since $\mathfrak q$ is unramified in $LK$, the order of the Frobenius automorphism $\mathrm{Frob}_{\mathfrak q}$ is equal to the inertia degree $f= f(\mathfrak Q/\mathfrak q)$ of $\mathfrak Q$ over $\mathfrak q$.  When the image $\rho(\mathrm{Gal}(LK/K))$ is isomorphic to $\mathrm{Gal}(LK/K)$, we know that the order of $\rho(\mathrm{Frob}_{\mathfrak q})$ is equal to the order of the $\mathrm{Frob}_{\mathfrak q}$.  Here the order of $\rho(\mathrm{Frob}_{\mathfrak q})$ is sufficient to determine the trace of $\rho(\mathrm{Frob}_{\mathfrak q})$.

Denote by $a_q$ the trace of $\rho_{\QQ}(\mathrm{Frob}_q)$.  When $q$ splits in $K$, say $q \cdot \mathcal O = \frak q \bar{\frak q}$, we have 
\[
	a_{\frak q} = a_{\bar{\frak q}} = a_q.
\]
Now consider primes $q \cdot \mathcal O = \frak q$ that are inert in $K$.  Suppose $a_q = \alpha + \beta$, for some $\alpha, \beta \in \bar \FF_{\ell}$ with $\alpha \beta = \det(\rho_{\QQ}(\mathrm{Frob}_q))$.  Then we have
\[
\begin{aligned}
	a_{\frak q} &= \alpha^2 + \beta^2 \\
				&= (\alpha + \beta)^2 - 2\det(\rho_{\QQ}(\mathrm{Frob}_q)) \\
				&= a_q^2 - 2\det(\rho_{\QQ}(\mathrm{Frob}_q)).
\end{aligned}
\]
Since the image of our representation is $\mathrm{SL}_2(\FF_3)$, we have $\det = 1$ and so we get 
\[
\begin{aligned}
	a_{\frak q} &= a_q^2 - 2 \\
		&\equiv a_q^2 + 1 \; \mod \; 3.
\end{aligned}
\]

In Table \ref{tab:reps} in Section \ref{tables} we give representatives for each conjugacy class in $\mathrm{SL}_2(\FF_3)$, the order of those elements, and $a_{\frak q} = \mathrm{tr}(\rho(\mathrm{Frob}_{\frak q}))$ for $q$ split and for $q$ inert in $K$.

To compute the weights, we look at the representation locally at $\ell = 3$ and apply the BDJ recipe.  We compute the ramification index $e_3 = 4$ and the inertia degree $f_3 = 2$.  The decomposition group at $3$ for $L$ over $\QQ$ has order $8$, and the only order $8$ subgroup is the quaternion group $Q_8$.  If we consider the restriction of $\rho_{\QQ}$ to the decomposition group, we would be in the irreducible case of the BDJ recipe.  However, we want to consider the restriction of $\rho_{\QQ}$ first to $K = \QQ(i)$, which we denote by $\rho$, and then to the decomposition group at $3$.  Then we have $\rho$ restricted to $D_3$ is reducible and, writing $\omega$ for a third root of unity in $\bar \FF_{\ell}^{\ast}$, we can write
\[
	\rho|_{D_3} \sim \twodmatrix{\omega^{(\ell^2 - 1)/4}}{0}{0}{\omega^{-(\ell^2 - 1)/4}}.
\]
The BDJ recipe includes a combinatorial analysis (see \cite[p. 18-24]{bdj} -- it is too long to reproduce here) that allows one to easily compute the predicted weights based on the power of $\omega$ in the above local description of $\rho$.  In Table \ref{tab:A4_61} we list these predicted weights.  The column labelled $B$ gives data relevant to this combinatorial analysis of in \cite{bdj}.  The $\vec a$ and $\vec b$ columns together determine one Serre weight.  In this case we have $\ell = 3$, which is inert in $K = \QQ(i)$, so we have only one prime $\frak p$ dividing $\ell$ but we have two embeddings $\tau_0, \tau_1 : k_{\frak p} \into \overline{\FF}_{\ell}$.  Letting $\vec a = (a_0, a_1)$ and $\vec b = (b_0, b_1)$ denote the corresponding exponents, our Serre weights are of the form
\begin{equation}\label{eq_Vp_inert}
V = V_{\mathfrak p} = \bigtensor_{i = 0, 1} \left ( \mathrm{det}^{a_i} \tensor_{k_{\mathfrak p}} \mathrm{Sym}^{b_i - 1} k_{\mathfrak p}^2 \right ) \tensor_{\tau_i} \bar \FF_\ell.
\end{equation}
The final column simply indicates that a corresponding modular form was indeed found for that weight.

\begin{table}[h]
\begin{center}
\caption{$A_4$ representation mod $3$ with level $\frak n = 61$}
\label{tab:A4_61}
\smallskip
\begin{tabular}{c|c|c|c|c|c}
	$\ell$ & level & $B$ & $\vec a = (a_0, a_1)$ & $\vec b = (b_0, b_1)$  & $f$ \\ \hline
	$3$ & $61$ & $\{ 0, 1 \}$ & $(0,2)$ & $(1,1)$ & $\checkmark$ \\
	$3$ & $61$ & $\{ 0, 1 \}$ & $(0,2)$ & $(3,3)$ & $\checkmark$ \\
	$3$ & $61$ & $\{ 1 \}$ & $(1,1)$ & $(2,2)$  & $\checkmark$ \\
	$3$ & $61$ & $\{ 0 \}$ & $(0,0)$ & $(2,2)$  & $\checkmark$ \\
	$3$ & $61$ & $\emptyset$ & $(2,0)$ & $(1,1)$  & $\checkmark$ \\
	$3$ & $61$ & $\emptyset$ & $(2,0)$ & $(3,3)$ & $\checkmark$ \\
\end{tabular}
\end{center}
\end{table}
In Table \ref{tab:examples} in Section \ref{tables} we list, for some small primes $\frak q$ of $\mathcal O$, the order of $\mathrm{Frob}_{\frak q}$ along with the coefficients $a_{\frak q}$ of the corresponding system of eigenvalues found.    From Table \ref{tab:reps}, we can determine $\mathrm{tr}(\rho(\mathrm{Frob}_{\frak q}))$ (from the order of $\mathrm{Frob}_{\frak q}$) and see that for each prime $\frak q$ this matches the eigenvalue found, taken mod $3$.

\subsection {Representations from Class Field Theory} \label{exsgalreps_cft}

The following examples are not arising as base change of even representations over $\QQ$, but are representations directly over $K=\QQ(i)$.  We get these examples by considering quadratic extensions of $\QQ(i)$ that are ramified only at a single prime $\frak p$, split over $\QQ$, and that have class group isomorphic to the cyclic group of order $3$.   In the following we will consider these representations mod $5$ and mod $7$.  This allows us to see the different behaviour in the weights modulo an inert prime as compared to a split prime.

The representation $\rho$ will factor through an extension $L$ of $K=\QQ(i)$ where $G = \mathrm{Gal}(L/K)$ is isomorphic to $D_3$.  There is one irreducible $2$-dimensional representation $\rho$ of the dihedral group $D_3$.  In Table \ref{tab:reps} in Section \ref{tables}, we list the images under $\rho$ of the elements of $D_3$, along with the order and trace of each.  In all cases the representation will have a quadratic character $\varepsilon : (\mathcal O/\frak p)^{\ast} \rightarrow \bar \FF_{\ell}^{\ast}$.  

We compute the level as described in Section \ref{exsgalreps_polys} above.  For these examples, we need only consider this single ramified prime $\frak p$, which will have ramification index $e_{\frak p} = 2$ in the extension $L$ over $K$.  The representation $\rho$ restricted to the order $2$ subgroup of $D_3$ fixes a one-dimensional subspace, so that the level in each case will be precisely equal to this prime, i.e., $\frak n = \frak p$.   The primes for which we found such dihedral extensions are 
\[
\begin{aligned}
	\frak n &= (8+17i) \mathrm{ \; (lying \; over \; p = 353)}, \\
	\frak n &= (13+28i) \mathrm{ \; (lying \; over \; p = 953) \; and } \\
	\frak n &= (8+35i) \mathrm{ \; (lying \; over \; p = 1289)}.  
\end{aligned}
\]

For these $D_3$ examples over $K = \QQ(i)$, we compute the values $\mathrm{tr}(\rho(\mathrm{Frob}_{\frak q}))$ for each prime $\frak q$ by computing the product of 
\begin{enumerate}
	\item the inertia degree of $\frak q \subset \mathcal O$ in the quadratic extension $L$, and 
	\item the order of a prime $\frak Q$ above $\frak q$ in the class group of $L$.  
\end{enumerate}
This product gives us the order of $\mathrm{Frob}_{\frak q}$ in the Galois group, which is isomorphic to $D_3$.  We denote this order by $o(\mathrm{Frob}_{\frak q})$.

In all three examples, the primes above $5$ and $7$ in $\QQ(i)$ will be unramified in the dihedral extension $L$ over $K$.  Thus, in all cases, the representation $\rho$ restricted to the inertia group $I$ at $\ell$ will be trivial.  The weight computations will depend on whether the prime $\ell$ is split or inert.  

In Table \ref{tab:D3overK_wts}, we list all the predicted weights for the above examples and indicate in the final column (under $f$) whether a corresponding modular form was found.   The weights do not depend on the level, only on the prime $\ell$.  For $\ell = 7$, which is inert in $K = \QQ(i)$, the $\vec a$ and $\vec b$ exponents correspond to the two embeddings of $k_{\frak p} \into \overline{\FF}_{\ell}$ as in Equation \eqref{eq_Vp_inert}.   For $\ell = 5$, which splits in $K = \QQ(i)$, we have two primes $\frak p_0$ and $\frak p_1$ dividing $\ell$.  The $\vec a = (a_0, a_1)$ and $\vec b = (b_0, b_1)$ then correspond to the two primes so that the Serre weights in this case are of the form
\[
V = \bigtensor_{i = 0, 1} \left ( \mathrm{det}^{a_i} \tensor_{k_{\mathfrak p_i}} \mathrm{Sym}^{b_i - 1} k_{\mathfrak p_i}^2 \right ) \tensor \bar \FF_\ell.
\]
\begin{table}[h!t]
\caption{$D_3$ representations mod $5$ and mod $7$ }
\label{tab:D3overK_wts}
\begin{center}
\begin{tabular}{c|c|c|c}
	$\ell$ & $\vec a = (a_0, a_1)$ & $\vec b = (b_0, b_1)$ & $f$\\ \hline
	$5$ & $(0,0)$ & $(4,4)$ & $\checkmark$ \\ \hline
	$7$ & $(0,0)$ & $(6,6)$ & $\checkmark$ \\
	$7$ & $(5,6)$ & $(1,7)$ & $\checkmark$ \\
	$7$ & $(6,5)$ & $(7,1)$ & $\checkmark$ \\ 
\end{tabular}
\end{center}
\end{table}

In Table \ref{tab:examples} in Section \ref{tables} we list, for some small primes $\frak q$, the order of $\mathrm{Frob}_{\frak q}$ (from which one can compute $\mathrm{tr}(\rho(\mathrm{Frob}_{\frak q}))$ as discussed above) along with the coefficients $a_{\frak q}$ of the systems of eigenvalues (considered mod $5$ and mod $7$).  These values are listed for each of the levels $\frak n = 8+17i$, $\frak n = 13+28i$ and $\frak n = 8+35i$.  In all cases, we found the corresponding systems of eigenvalues in all weights predicted for both $\ell = 5$ and $\ell = 7$.

\section {Computations} \label{tables}

Table \ref{tab:reps} gives representatives of the groups $\mathrm{SL}_2(\FF_3)$ and of $D_3$ (for the $D_3$ elements,  $\omega $ is a third root of unity in $\bar \FF_{\ell}^{\ast}$) along with the order and trace of each element.  In the case of $\mathrm{SL}_2(\FF_3)$, we get the representation by taking the base change to $K = \QQ(i)$ of an even representation over $\QQ$ and so we give the trace for $q$ split and for $q$ inert (in $K$) as described above in Section \ref{exsgalreps_polys}.
\begin{table}[h!t]
\caption{Representatives of $\mathrm{SL}_2(\FF_3)$ and of $D_3$}
\label{tab:reps}
\begin{center}
\begin{tabular}{|c|c|c|c||c|c|c|} \hline
	\multicolumn{4}{|c||}{$\mathrm{SL}_2(\FF_3)$} 	&	\multicolumn{3}{c|}{$D_3$} 	\\ \hline
	\multirow{2}{*}{Rep} & \multirow{2}{*}{Order} & Trace & Trace & \multirow{2}{*}{Rep}		& \multirow{2}{*}{Order} 	& \multirow{2}{*}{Trace} \\ 
	 &  & ($q$ split) & ($q$ inert) 					&		&  	&  \\ \hline
	$\twodmatrix{1}{0}{0}{1}$ & $1$ & $-1$ & $-1$ 		&	$\twodmatrix{1}{0}{0}{1}$	& $1$	& $2$	\\
	$\twodmatrix{-1}{0}{0}{-1}$ & $2$ & $1$ & $-1$ 	&	$\twodmatrix{\omega}{0}{0}{\omega ^{-1}}$	& $3$	& $-1$	\\
	$\twodmatrix{0}{1}{-1}{-1}$ & $3$ & $-1$ & $-1$ 	&	$\twodmatrix{\omega ^2}{0}{0}{\omega ^{-2}}$	& $3$	& $-1$ 	\\
	$\twodmatrix{0}{-1}{1}{-1}$ & $3$ & $-1$ & $-1$ 	&	$\twodmatrix{0}{1}{1}{0}$			& $2$	& $0$	\\
	$\twodmatrix{0}{-1}{1}{0}$ & $4$ & $0$ & $1$ 		&	$\twodmatrix{0}{\omega ^{-1}}{\omega}{0}$	& $2$	& $0$ 	\\
	$\twodmatrix{0}{-1}{1}{1}$ & $6$ & $1$ & $-1$ 		&	$\twodmatrix{0}{\omega ^{-2}}{\omega ^2}{0}$	& $2$	& $0$ 	\\
	$\twodmatrix{0}{1}{-1}{1}$ & $6$ & $1$ & $-1$ 		&	&	&	\\ \hline
\end{tabular}
\end{center}
\end{table}

For each of the representations computed in Section \ref{galreps} above, we have computed the orders of Frobenius elements for some small primes $\frak q$.  These are listed in Table \ref{tab:examples}.  From the orders of Frobenius elements, one can find $\mathrm{tr}(\rho(\mathrm{Frob}_{\frak q}))$ from Table \ref{tab:reps}.  Table \ref{tab:examples} also lists the eigenvalues $a_{\frak q}$ of modular forms found using the author's code.  One can then check that the traces match up with the eigenvalues modulo the appropriate primes (listed in the headers of the tables).  The modular forms were found in all of the predicted weights, as listed in Tables \ref{tab:D3overK_wts} and \ref{tab:A4_61} in Section \ref{galreps}.

\begin{table}[h!t]
\centering
{
\caption{Examples of Galois Representations and Corresponding Modular Forms}
\bigskip
\label{tab:examples}
\begin{tabular}{|c||c|c||c|c||c|c||c|c|} \hline
		& \multicolumn{2}{c||}{$A_4$, $\ell = 3$} & \multicolumn{2}{c||}{$D_3$, $\ell = 5,7$} & \multicolumn{2}{c||}{$D_3$, $\ell = 5,7$} & \multicolumn{2}{c|}{$D_3$, $\ell = 5,7$} \\  \hline
		& \multicolumn{2}{c||}{$\frak n = 61$} & \multicolumn{2}{c||}{$\frak n = 8+17i$} & \multicolumn{2}{c||}{$\frak n = 13+28i$} & \multicolumn{2}{c|}{$\frak n = 8+35i$} \\  \hline
	$\frak q$ & $o(\mathrm{Frob}_{\frak q})$ & $a_{\frak q}$ & $o(\mathrm{Frob}_{\frak q})$ & $a_{\frak q}$ & $o(\mathrm{Frob}_{\frak q})$ & $a_{\frak q}$ & $o(\mathrm{Frob}_{\frak q})$ & $a_{\frak q}$  \\ \hline 
	$1+i$ & $4$ & $0$	& $3$ & $-1$ & $3$ & $-1$ & $2$ & $0$ \\
	$1-2i$ & $6$ & $1$  & $2$ & $0$ & $2$ & $0$  & $3$ & $-1$ \\
	$1+2i$ & $6$ & $1$  & $3$ & $-1$ &  $3$ & $-1$ & $3$ & $-1$ \\ 
	$3$ & - & - & $2$ & $0$ & $2$ & $0$ & $2$ & $0$   \\
	$2-3i$ & $6$ & $1$  & $3$ & $-1$ & $1$ & $2$ & $2$ & $0$  \\
	$2+3i$ & $6$ & $1$  & $2$ & $0$ & $3$ & $-1$ & $3$ & $-1$ \\
	$1-4i$ & $6$ & $1$  & $3$ & $-1$ & $2$ & $0$ & $3$ & $-1$ \\
	$1+4i$ & $6$ & $1$  & $3$ & $-1$ &  $2$ & $0$ & $2$ & $0$  \\
	$2-5i$ & $3$ & $2$  & $2$ & $0$ & $3$ & $-1$ & $2$ & $0$  \\
	$2+5i$ & $3$ & $2$ & $2$ & $0$ & $3$ & $-1$ & $2$ & $0$  \\
	$1-6i$ & $4$ & $0$  & $3$ & $-1$ & $3$ & $-1$ & $3$ & $-1$  \\ 
	$1+6i$ & $4$ & $0$  & $2$ & $0$ &  $3$ & $-1$ & $2$ & $0$ \\
	$4-5i$ & $4$ & $0$  & $2$ & $0$ & $2$ & $0$ & $3$ & $-1$   \\
	$4+5i$ & $4$ & $0$  & $2$ & $0$ & $2$ & $0$ & $1$ & $2$   \\
	$7$ & $3$ & $2$  	& $2$ & $0$ &  $3$ & $-1$ & $3$ & $-1$   \\
	$2-7i$ & $4$ & $0$   & $2$ & $0$ &  $2$ & $0$ & $3$ & $-1$ \\
	$2+7i$ & $4$ & $0$  & $1$ & $2$ &  $2$ & $0$ & $1$ & $2$ \\
	$5-6i$ & - & - & $2$ & $0$ & $3$ & $-1$ & $2$ & $0$ \\
	$5+6i$ & - & - & $2$ & $0$ & $2$ & $0$ & $3$ & $-1$ \\ 
	$3-8i$ & $3$ & $2$ & $2$ & $0$ & $2$ & $0$ & $3$ & $-1$ \\
	$3+8i$ & $3$ & $2$ & $2$ & $0$ &  $2$ & $0$ & $1$ & $2$ \\
	$5-8i$ & $4$ & $0$ & $2$ & $0$ &  $2$ & $0$ & $2$ & $0$ \\
	$5+8i$ & $4$ & $0$ & $1$ & $2$ &  $3$ & $-1$  & $3$ & $-1$ \\
	$4-9i$ & $3$ & $2$ & $3$ & $-1$ &  $3$ & $-1$ & $3$ & $-1$  \\
	$4+9i$ & $3$ & $2$ & $1$ & $2$   &  $2$ & $0$ & $2$ & $0$  \\
	$1-10i$ & $6$ & $1$ & $3$ & $-1$ &  $1$ & $2$ & $1$ & $2$   \\
	$1+10i$ & $6$ & $1$ & $2$ & $0$ &  $2$ & $0$ & $3$ & $-1$ \\
	$3-10i$ & $6$ & $1$ & $3$ & $-1$ &  $2$ & $0$ & $2$ & $0$ \\
	$3+10i$ & $6$ & $1$ & $3$ & $-1$ &  $2$ & $0$ & $3$ & $-1$ \\
	$7-8i$ & $4$ & $0$ & $2$ & $0$ &  $2$ & $0$ & $1$ & $2$ \\
	$7+8i$ & $4$ & $0$ & $2$ & $0$ &  $2$ & $0$ & $2$ & $0$ \\
	$11$ & $4$ & $1$ & $3$ & $-1$ &  $2$ & $0$  & $2$ & $0$ \\
	$4-11i$ & $6$ & $1$ & $3$ & $-1$ &  $3$ & $-1$ & $2$ & $0$ \\
	$4+11i$ & $6$ & $1$ & $2$ & $0$ &  $2$ & $0$ & $2$ & $0$ \\
	$7-10i$ & $4$ & $0$ & $2$ & $0$ &  $2$ & $0$ & $2$ & $0$ \\
	$7+10i$ & $4$ & $0$ & $3$ & $-1$ &  $1$ & $2$  & $3$ & $-1$ \\
	\hline
\end{tabular}
}
\end{table}

 \pagebreak

\section {Final Remarks} \label{final_rmks}

In summary, the computational evidence presented here supports the surmise that the BDJ conjectural weight recipe for totally real fields will hold in the case of imaginary quadratic fields as well.  For the examples of Galois representations computed here, corresponding modular forms were found in all of the predicted weights.   Of the examples computed in \cite{torrey}, there were two exceptions. In one example we did not find the form in all the weights because the computations were too large for the program, so we have not yet looked for all of the predicted weights.  In the other example, we found the form in only two of the four predicted weights.  It does look like we may have found a twist of the form in the remaining two weights.  It is not yet understood why only two of those four weights were found.

The modular symbols computation method used in my program is justified here only for $K = \QQ(i)$.  I expect it will be straightforward to use the same methodology for all the other Euclidean class number one imaginary quadratic fields.  Note, however, that for each field one must justify an algebraic proposition such as Proposition \ref{prop_algprop} presented in Section \ref{alg_prop}.  We already know the relations we expect to work for each of these fields, namely those relations computed by Cremona, et al.  See, for example, \cite{cremonaB}. 

Several students of Cremona have extended the modular symbols method (for trivial weights) to imaginary quadratic fields of higher class number.  I expect, with some work, that the method presented in this thesis can be joined with their work to compute modular forms with arbitrary weight for imaginary quadratic fields of higher class number.

\section {Acknowledgements} \label{acknowledgements}

The work presented in this paper is from my thesis, written at King's College London under the supervision of Fred Diamond.  I would like to thank him for all his help.  I would also like to gratefully acknowledge helpful discussions with Avner Ash, John Cremona and Gabor Wiese.





\bibliographystyle{model1b-num-names}
\bibliography{numtheory}

\begin{thebibliography}{28}
\expandafter\ifx\csname natexlab\endcsname\relax\def\natexlab#1{#1}\fi
\providecommand{\bibinfo}[2]{#2}
\ifx\xfnm\relax \def\xfnm[#1]{\unskip,\space#1}\fi
\bibitem[{Ash(1994)}]{ash}
\bibinfo{author}{A.~Ash}, \bibinfo{title}{Unstable cohomology of
  $\text{SL}(n,\mathcal{O})$}, \bibinfo{journal}{J. of Algebra}
  \bibinfo{volume}{167} (\bibinfo{year}{1994}) \bibinfo{pages}{330--342}.
\bibitem[{Ash et~al.(2002)Ash, Doud and Pollack}]{ashDP}
\bibinfo{author}{A.~Ash}, \bibinfo{author}{D.~Doud},
  \bibinfo{author}{D.~Pollack}, \bibinfo{title}{{G}alois representations with
  conjectural connections to arithmetic cohomology}, \bibinfo{journal}{Duke
  Math. J.} \bibinfo{volume}{112} (\bibinfo{year}{2002})
  \bibinfo{pages}{521--579}.
\bibitem[{Ash et~al.(2011)Ash, Gunnells and McConnell}]{agm}
\bibinfo{author}{A.~Ash}, \bibinfo{author}{P.~Gunnells},
  \bibinfo{author}{M.~McConnell}, \bibinfo{title}{Torsion in the cohomology of
  congruence subgroups of $\text{SL}(4, z)$ and {G}alois representations},
  \bibinfo{journal}{J. of Algebra} \bibinfo{volume}{325} (\bibinfo{year}{2011})
  \bibinfo{pages}{404--415}.
\bibitem[{Borel and Serre(1973)}]{borelserre}
\bibinfo{author}{A.~Borel}, \bibinfo{author}{J.P. Serre},
  \bibinfo{title}{Corners and arithmetic groups}, \bibinfo{journal}{Comment.
  Math. Helv.} \bibinfo{volume}{48} (\bibinfo{year}{1973})
  \bibinfo{pages}{436--491}.
\bibitem[{Bosma et~al.(1997)Bosma, Cannon and Playoust}]{magma}
\bibinfo{author}{W.~Bosma}, \bibinfo{author}{J.~Cannon},
  \bibinfo{author}{C.~Playoust}, \bibinfo{title}{The {M}agma algebra system.
  {I}. {T}he user language}, \bibinfo{journal}{J. of Symb. Comput.}
  \bibinfo{volume}{24} (\bibinfo{year}{1997}) \bibinfo{pages}{235--265}.
  \bibinfo{note}{Available from {http://magma.maths.usyd.edu.au/magma/} as of
  June 24, 2011}.
\bibitem[{Buzzard et~al.(ress)Buzzard, Diamond and Jarvis}]{bdj}
\bibinfo{author}{K.~Buzzard}, \bibinfo{author}{F.~Diamond},
  \bibinfo{author}{F.~Jarvis}, \bibinfo{title}{On {S}erre's conjecture for mod
  $\ell$ {G}alois representations over totally real fields},
  \bibinfo{journal}{Duke Math. J.}  (\bibinfo{year}{in press}).
\bibitem[{Cremona(1984)}]{cremonaB}
\bibinfo{author}{J.~Cremona}, \bibinfo{title}{Hyperbolic tessellations, modular
  symbols, and elliptic curves over complex quadratic fields},
  \bibinfo{journal}{Compos. Math.} \bibinfo{volume}{51} (\bibinfo{year}{1984})
  \bibinfo{pages}{275--323}.
\bibitem[{\c{S}eng\"{u}n(2008)}]{sengun}
\bibinfo{author}{M.H. \c{S}eng\"{u}n}, \bibinfo{title}{Serre's Conjecture over
  Imaginary Quadratic Fields}, Ph.D. thesis, University of Wisconsin-Madison,
  \bibinfo{year}{2008}.
\bibitem[{Demb\'{e}l\'{e} et~al.(ults)Demb\'{e}l\'{e}, Diamond and
  Roberts}]{ddr}
\bibinfo{author}{L.~Demb\'{e}l\'{e}}, \bibinfo{author}{F.~Diamond},
  \bibinfo{author}{D.~Roberts}, \bibinfo{title}{Numerical examples and evidence
  for {S}erre's conjecture over totally real fields},
  \bibinfo{year}{unpublished results}.
\bibitem[{Dieulefait(2004)}]{dieulefait}
\bibinfo{author}{L.~Dieulefait}, \bibinfo{title}{Existence of families of
  {G}alois representations and new cases of the {F}ontaine-{M}azur conjecture},
  \bibinfo{journal}{J. Reine Angew. Math.} \bibinfo{volume}{577}
  (\bibinfo{year}{2004}) \bibinfo{pages}{147--151}.
\bibitem[{Figueiredo(1999)}]{figueiredo}
\bibinfo{author}{L.M. Figueiredo}, \bibinfo{title}{Serre's conjecture for
  imaginary quadratic fields}, \bibinfo{journal}{Compos. Math.}
  \bibinfo{volume}{118} (\bibinfo{year}{1999}) \bibinfo{pages}{103--122}.
\bibitem[{Gee(2007)}]{gee1}
\bibinfo{author}{T.~Gee}, \bibinfo{title}{Companion forms over totally real
  fields {II}}, \bibinfo{journal}{Duke Math. J.} \bibinfo{volume}{136}
  (\bibinfo{year}{2007}) \bibinfo{pages}{275--284}.
\bibitem[{Gee(2011)}]{gee2}
\bibinfo{author}{T.~Gee}, \bibinfo{title}{On the weights of mod p {H}ilbert
  modular forms}, \bibinfo{journal}{Invent. Math,} \bibinfo{volume}{184}
  (\bibinfo{year}{2011}) \bibinfo{pages}{1--46}.
\bibitem[{Khare and Wintenberger(2009{\natexlab{a}})}]{khare_wintenberger_I}
\bibinfo{author}{C.~Khare}, \bibinfo{author}{J.P. Wintenberger},
  \bibinfo{title}{{S}erre's modularity conjecture ({I})},
  \bibinfo{journal}{Invent. Math.} \bibinfo{volume}{178}
  (\bibinfo{year}{2009}{\natexlab{a}}) \bibinfo{pages}{485--504}.
\bibitem[{Khare and Wintenberger(2009{\natexlab{b}})}]{khare_wintenberger_II}
\bibinfo{author}{C.~Khare}, \bibinfo{author}{J.P. Wintenberger},
  \bibinfo{title}{{S}erre's modularity conjecture ({II})},
  \bibinfo{journal}{Invent. Math.} \bibinfo{volume}{178}
  (\bibinfo{year}{2009}{\natexlab{b}}) \bibinfo{pages}{505--586}.
\bibitem[{Kisin(2009{\natexlab{a}})}]{kisin1}
\bibinfo{author}{M.~Kisin}, \bibinfo{title}{Modularity of $2$-adic
  {B}arsotti-{T}ate representations}, \bibinfo{journal}{Invent. Math.}
  \bibinfo{volume}{178} (\bibinfo{year}{2009}{\natexlab{a}})
  \bibinfo{pages}{587--634}.
\bibitem[{Kisin(2009{\natexlab{b}})}]{kisin2}
\bibinfo{author}{M.~Kisin}, \bibinfo{title}{Moduli of finite flat group schemes
  and modularity}, \bibinfo{journal}{Ann. of Math.} \bibinfo{volume}{170}
  (\bibinfo{year}{2009}{\natexlab{b}}) \bibinfo{pages}{1085--1180}.
\bibitem[{Kl\"{u}ners and Malle(site)}]{db_km}
\bibinfo{author}{J.~Kl\"{u}ners}, \bibinfo{author}{G.~Malle}, \bibinfo{title}{A
  database for number fields}, \bibinfo{year}{website}.
  \bibinfo{note}{Available from
  {http://www.math.uni-duesseldorf.de/$\sim$klueners/minimum/minimum.html} as
  of June 24, 2011}.
\bibitem[{Martin(2001)}]{martin}
\bibinfo{author}{F.~Martin}, \bibinfo{title}{P\'{e}riodes de formes modulaires
  de poids 1}, Ph.D. thesis, Universit\'{e} Paris 7, \bibinfo{year}{2001}.
\bibitem[{PARI(2005)}]{PARI2}
PARI, \bibinfo{title}{{PARI/GP, version {\tt 2.3.1}}}, \bibinfo{year}{2005}.
  \bibinfo{note}{Available from {http://pari.math.u-bordeaux.fr/} as of June
  24, 2011}.
\bibitem[{Serre(1987)}]{serre2}
\bibinfo{author}{J.P. Serre}, \bibinfo{title}{Sur les repr\'{e}sentations
  modulaires de degr\'{e} 2 de $\text{Gal}(\bar {\QQ}/ {\QQ})$},
  \bibinfo{journal}{Duke Math. J.} \bibinfo{volume}{54} (\bibinfo{year}{1987})
  \bibinfo{pages}{179--230}.
\bibitem[{Stein(2007)}]{stein}
\bibinfo{author}{W.~Stein}, \bibinfo{title}{Modular Forms, a Computational
  Approach}, volume~\bibinfo{volume}{79} of \textit{\bibinfo{series}{Graduate
  Studies in Mathematics}}, \bibinfo{publisher}{American Mathematical Society},
  \bibinfo{address}{Providence}, \bibinfo{year}{2007}.
\bibitem[{Stein(site)}]{sage}
\bibinfo{author}{W.~Stein}, \bibinfo{title}{Sage mathematical software system},
  \bibinfo{year}{website}. \bibinfo{note}{Available from
  {http://www.sagemath.org/} as of June 24, 2011}.
\bibitem[{Taylor(2002)}]{taylor}
\bibinfo{author}{R.~Taylor}, \bibinfo{title}{Remarks on a conjecture of
  {F}ontaine and {M}azur}, \bibinfo{journal}{Inst. Math. Jussieu}
  \bibinfo{volume}{1} (\bibinfo{year}{2002}) \bibinfo{pages}{1--19}.
\bibitem[{Taylor and Wiles(1995)}]{taylor_wiles}
\bibinfo{author}{R.~Taylor}, \bibinfo{author}{A.~Wiles},
  \bibinfo{title}{Ring-theoretic properties of certain {H}ecke algebras},
  \bibinfo{journal}{Ann. of Math.} \bibinfo{volume}{141} (\bibinfo{year}{1995})
  \bibinfo{pages}{553--572}.
\bibitem[{Torrey(2009)}]{torrey}
\bibinfo{author}{R.~Torrey}, \bibinfo{title}{On {S}erre's Conjecture over
  Imaginary Quadratic Fields}, Ph.D. thesis, King's College London,
  \bibinfo{year}{2009}.
\bibitem[{Wiese(2009)}]{wiese_hecke}
\bibinfo{author}{G.~Wiese}, \bibinfo{title}{On modular symbols and the
  cohomology of {H}ecke triangle surfaces}, \bibinfo{journal}{International J.
  of Number Theory} \bibinfo{volume}{5} (\bibinfo{year}{2009})
  \bibinfo{pages}{89--108}.
\bibitem[{Wiles(1995)}]{wiles}
\bibinfo{author}{A.~Wiles}, \bibinfo{title}{Modular elliptic curves and
  {F}ermat's {L}ast {T}heorem}, \bibinfo{journal}{Ann. of Math.}
  \bibinfo{volume}{141} (\bibinfo{year}{1995}) \bibinfo{pages}{443--551}.

\end{thebibliography}







\end{document}